\documentclass[11pt,reqno,oneside]{article}
\usepackage{amsmath}
\usepackage{amssymb}
\usepackage{url}
\usepackage{natbib}
\usepackage[dvips]{graphicx}
\usepackage[font=small,labelfont=bf]{caption}
\usepackage{pstricks}

\newtheorem{prop}{Proposition}

\title{Exact enumeration of cherries and pitchforks in ranked trees under the coalescent model}
\author{Filippo Disanto$^*$\and Thomas Wiehe\thanks{Institut f\"ur Genetik, Universit\"at zu K\"oln; Z\"ulpicher Stra\ss e 47a, 50674 K\"oln, Germany}}
\date{\today}

\begin{document}

\maketitle

\begin{abstract}
We consider exact enumerations and probabilistic properties of \emph{ranked} trees when generated under the random coalescent process. Using a new approach (see \cite{Steel99,rosenberg}), based on generating functions, we derive several statistics such as the exact probability of finding $k$ \emph{cherries} in a ranked tree of fixed size $n$. We then extend our method to consider also the number of \emph{pitchforks}. We find a recursive formula to calculate the joint and conditional probabilities of cherries and pitchforks when the size of the tree is fixed.
\end{abstract}

\section{Introduction}

Given a direction by time, ancestry relationship between species, individuals, alleles or cells can be depicted as a rooted tree. Of particular interest are binary rooted unordered trees. These can be further classified into several subclasses. Here we will 
%\emph{shape trees}, \emph{phylogenetic trees}, 
\emph{ranked trees}, which are defined below. 
%and \emph{labelled ranked trees}.

%In this work we focus on several combinatorial aspects concerning these classes of trees. First, we present their enumerations. If trees are generated under a random process it is of interest to know their probability distributions. 

We assume that trees are generated by the coalescent process. 
%and we derive the probabilities of several tree properties. 

An important parameter is the number of \emph{cherries} of a tree. By a new approach based on generating functions  we extend previous results (see for example \cite{Steel99}) deriving an exact formula for the probability of finding $k$ cherries in a ranked tree of size $n$. Furthermore, we show that several known statistics (see \cite{rosenberg}) concerning \emph{pitchforks} follow as corollaries from a partial differential equation which also gives an efficient recursion to compute the conditional probability distribution of pitchforks given a certain number of cherries.

%with the aim of answering the question \emph{how does a 'typical' tree look like?}

One motivation for this study comes from population genetics and the question how 'typical' \emph{coalescent trees} \cite{Wakeleybook} look like. Our results give some insight into structural properties of trees generated under the standard neutral model \cite{Tajima6628982}. These results provide a reference against which non neutral and/or non independently generated trees may be compared. To illustrate the latter we pay attention to trees which are linked along a recombining chromosome.

%recombination between homologous chromosomes (e.g. \cite{Hudson90}) affects the shape of the ancestral tree of a sample of $n$ genes. To compare the effect of recombination on linked trees to completely unlinked trees we will also provide a formula which gives the background probability of generating the same ranked tree as a result of two independent random coalescent processes.

\section{Preliminaries}

\begin{figure}
\begin{center}
\begin{tabular}{|c|c|c|}\hline
trees & \# cherries & \# pitchforks \\ \hline
\includegraphics*[angle=0,scale=.33,trim=0 0 0 0]{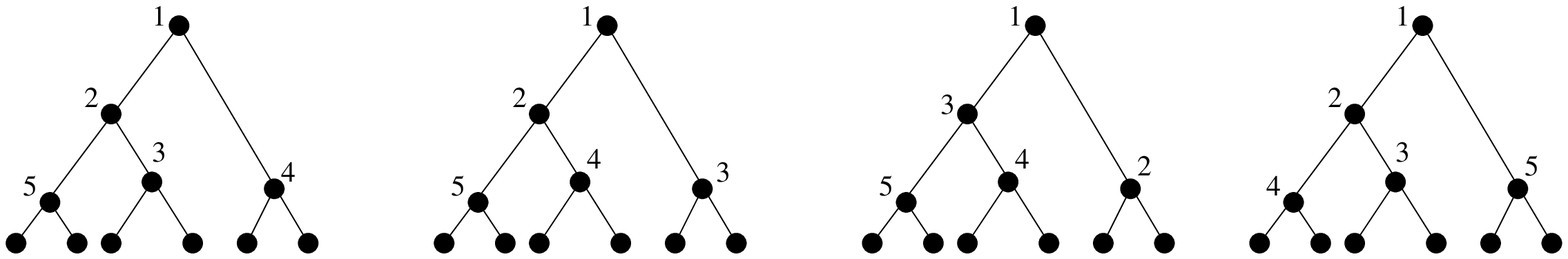} & 3 & 0 \\ \hline
\includegraphics*[angle=0,scale=.33,trim=0 0 0 0]{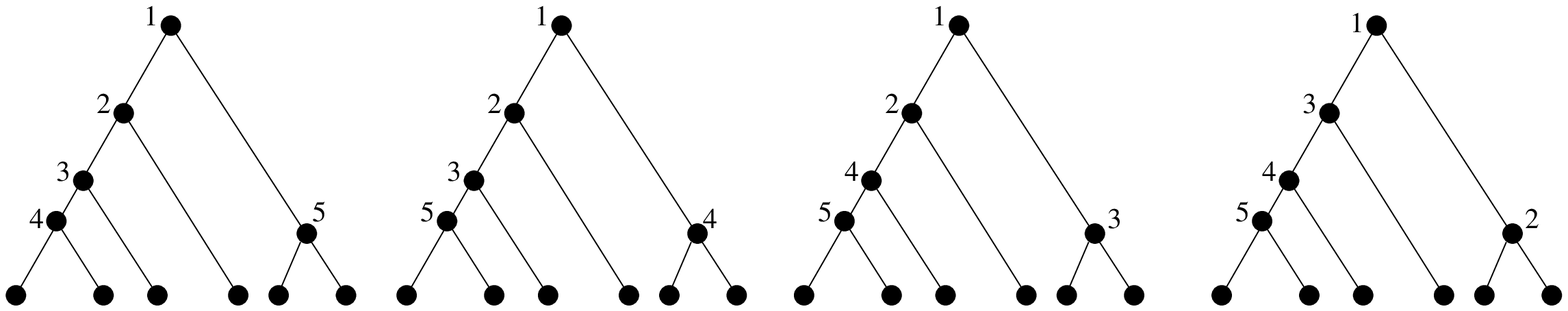} & 2 & 1 \\ \hline
\includegraphics*[angle=0,scale=.33,trim=0 0 0 0]{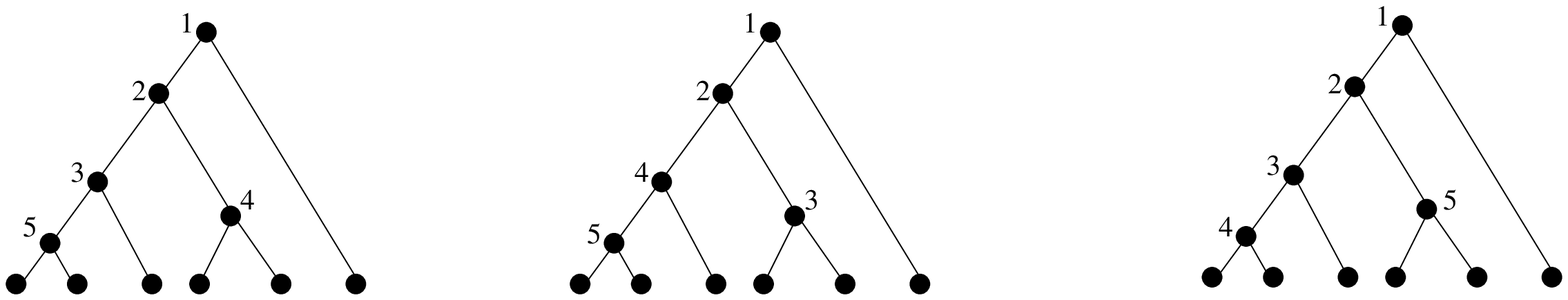} & 2 & 1 \\ \hline
\includegraphics*[angle=0,scale=.33,trim=0 0 0 0]{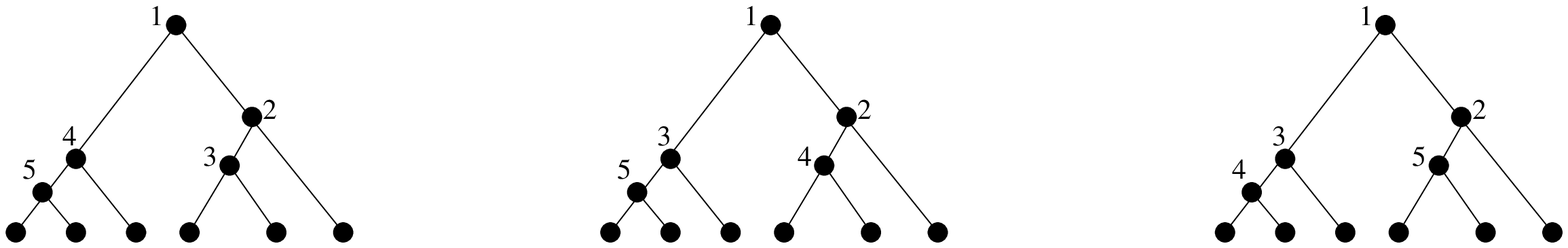} & 2 & 2 \\ \hline
\includegraphics*[angle=0,scale=.33,trim=0 0 0 0]{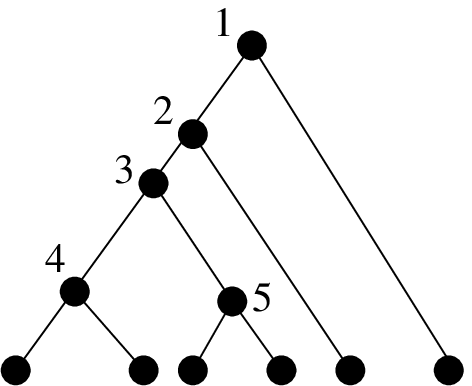} & 2 & 0 \\ \hline
\includegraphics*[angle=0,scale=.33,trim=0 0 0 0]{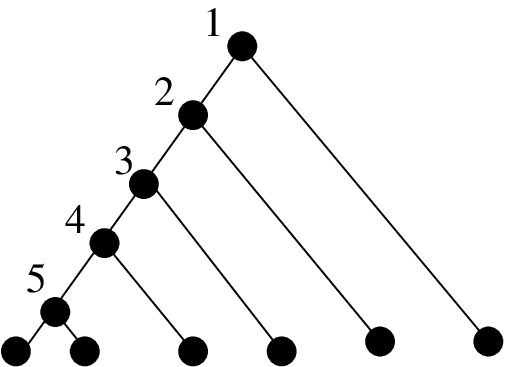} & 1 & 1 \\ \hline
\end{tabular}
\end{center}
\caption{The sixteen possible ranked trees of size six classified by shape. Within each class all possible orderings of the internal nodes are displayed. Number of cherries and pitchforks are indicated.}\label{ranked6}
\end{figure}

We start with some basic definitions.
A {\it binary rooted} tree is a tree with a root and in which all nodes have outdegree either $0$ or $2$.
Nodes with outdegree $2$ are called {\it internal}, nodes with outdegree $0$ are {\it external}. External nodes are also called {\it leaves}.
The size $n$ of a tree is the number of its external nodes.
The {\it subtree} of an internal node $i$ is the tree with root $i$. 
A tree is said to be \emph{un-ordered} when it is taken in the graph theoretic sense so that subtrees stemming from an internal node have not a left-right order between themselves. 
%The {\it left} ({\it right}) subtree of $i$
%is the subtree with a root which is the left (right) child node of $i$. In an {\it ordered} ({\it un-ordered}) rooted tree the left and right subtrees of an internal node are (not) ordered.
Here, we care about tree topology and we do not care about branch lengths. 
We consider the following class. 
A binary un-ordered tree of size $n$ is said to be a \emph{ranked tree} if  the set of internal nodes is totally ordered by labels belonging to $\{1,2,...,n \}$ in such a way that each child's label is greater than its parent's label, (see Fig.~\ref{ranked6}). The total order of internal labels can be interpreted as a historical time order; accordingly,  \citet{Harding1971} calls such trees \emph{histories}.

%\begin{itemize}
%\item[(i)] unlabelled trees; we will refer to these trees as \emph{shape trees}, see Fig.~\ref{examples}~(a). The name, for instance used in an early paper by \citet{Harding1971}, reflects the fact that the only information is its shape
%\item[(ii)] trees with totally ordered labels at internal nodes; more precisely, internal nodes are labelled in such a way that each child's label is greater than its parent's label; we will refer to these trees as \emph{ranked trees}, see Fig.~\ref{examples}~(b). The total order of internal labels can be interpreted as a historical time order; accordingly,  \citet{Harding1971} calls such trees \emph{histories}
%\item[(iii)] trees with labelled leaves; we will refer to these trees as \emph{phylogenetic trees} as done for example in \cite{Bona:1154505}, see Fig.~\ref{examples}~(c). Labels may be species names, for instance. \citet{Harding1971} calls this class  \emph{labelled shape trees}
%\item[(iv)] trees with external and ordered internal labels; we will refer to these trees as \emph{labelled ranked trees}, see Fig.~\ref{examples}~(d). \citet{Harding1971} calls them \emph{labelled histories}.
%\end{itemize}

%\begin{figure}
%\begin{center}
%\includegraphics*[scale=.40,trim=0 0 0 0]{esempialberi.eps}
%\end{center}
%\caption{(a) Shape tree, (b) ranked tree, (c) phylogenetic tree, (d) labelled ranked tree; all trees have size $n=6$.}\label{examples}
%\end{figure}

%\begin{array}{|c|c|c|}\hline
%Trees & cherries & pitchforks \hline
%\includegraphics*[angle=0,scale=.33,trim=0 0 0 0]{riga1.eps} & 3 & 0
%\hline
%\end{array}

%\begin{figure}
%\begin{center}
%\includegraphics*[scale=.4,trim=0 0 0 0]{alberi2.eps}
%\end{center}
%\caption{The sixteen ranked trees of size $6$. Each row contains the trees which have identical shapes. There are six different shapes for trees of size $n=6$. Note, there is one shape with 3 cherries (top row) and one shape with 1 cherry (bottom row). The remaining ones have 2 cherries.}\label{ranked6}
%\end{figure}

%We will denote by $\mathcal{S}$ (resp. $\mathcal{P},\mathcal{R},\mathcal{L}$) the set of shape (resp. phylogenetic, ranked, labelled ranked) trees and by $\mathcal{S}_n$ (resp. $\mathcal{P}_n,\mathcal{R}_n,\mathcal{L}_n$) the set of trees of size $a$.

We will denote by $\mathcal{R}$ the set of ranked trees and by $\mathcal{R}_n$ the set of trees of size $n$. In what follows, $n=n(t)$ always represents the number of leaves of a ranked tree $t$.

\bigskip
 
The cardinality of the set $\mathcal{R}_n$  is given by the following exponential generating function 

\begin{equation}
\mathcal{R}(x)  =  \sum_{n \geq 0} \frac{|\mathcal{R}_n|}{n !} x^n = \sec (x)+ \tan (x). \label{enumR}
\end{equation}

%\begin{eqnarray}
% \mathcal{S} (x) & = & x+\frac{1}{2}(\mathcal{S} (x)^2+\mathcal{S} (x^2)) \label{enumS}\\
% \mathcal{P}_n & = & \prod_{i=2}^{n}(2i-3) \label{enumP}\\
% \mathcal{R}(x) & = & \sec (x)+ \tan (x) \label{enumR}\\
% \mathcal{L}_n & = & \frac{n!(n-1)!}{2^{n-1}} \label{enumL}
%\end{eqnarray}

%In eq (\ref{enumS}) the function $\mathcal{S}(x)$ is the ordinary generating function of shape trees. Solving iteratively the equation one finds that the first terms of the sequence are (with $n>0$): 
%$$
%1,1,1,2,3,6,11,23,46,...
%$$ 
%and they correspond to sequence $A001190$ in \citet{MR1992789}. These numbers are known as \emph{Wedderburn-Etherington numbers}.

%In eq (\ref{enumP}) the first terms of the sequence are (with $n>0$): 
%$$
%1,1,3,15,105,945,...
%$$ 
%and they correspond to sequence $A001147$ in \citet{MR1992789}.

whose first coefficients  $|\mathcal{R}_n|$ (with $n>0$) are 
$$
1,1,1,2,5,16,61,272,... .
$$ 
Ranked trees can be bijectively mapped to $0$-$1$-$2$-\emph{increasing trees} (see Callan, 2005; \url{http://www.stat.wisc.edu/~callan/notes}). From this, it follows that the numbers  given by (\ref{enumR}) correspond to sequence $A000111$ in \citet{MR1992789} and are known as \emph{Euler numbers}.

%Finally, in eq (\ref{enumL}) the first terms of the sequence are (with $n>0$): 
%$$
%1,1,3,18,180,2700,...
%$$ 
%and they correspond to sequence $A006472$ in \cite{MR1992789}.

\subsection{Trees as a result of the coalescent process} \label{process}

The coalescent of size $n$ is a model for the genealogical history of a sample of $n$ genes. It has been introduced in population genetics by Kingman and Ewens \cite{kingman:1982,Kingman11102348} and has nowadays textbook status   \cite{Wakeleybook}. Ranked trees can be generated by the \emph{coalescent process}, which starts with $n$ leaves and works by successively coalescing two randomly chosen branches until it reaches the 'most recent common ancestor' when the last two remaining branches are joined. 
%The timing of coalescent events is governed by an exponential process. The genealogical history of a sample of homologous genes can be considered as the result of such a random process.
%It generates a ranked tree starting from the leaves and goes 'backward' in time. Each coaelscent event creates a new internal node of the tree. After completion there are $n-1$ internal nodes.
%\begin{figure}
%\begin{center}
%\includegraphics*[scale=.37,trim=0 0 0 0]{coalesce.eps}
%\end{center}
%\caption{A ranked tree of size $5$ as the result of a coalescent process.}\label{coalescent}
%\end{figure}

To reflect time order one can assign an integer to each internal node when created, for instance the label $n-1$ to the first coalescent event and $1$ to the last event, the most recent common ancestor, or the root of the tree. 
%(see Fig.~\ref{coalescent}).
%Note, that a ranked tree can also be obtained as the result of a random process moving forward in time. 
%The forward analogon to the colaescent process is the pure birth process or {\it Yule} process.
%In the next section we ask for the probabilities of trees generated under the coalescent process.

%\subsection{\small{Probability distribution of trees}}\label{weights}
The probability distribution of ranked trees $P_{\mathcal{R}}$ generated under the coalescent process is essentially contained in the paper of \citet{Tajima6628982} and it is described below. 

%Removing labels from a ranked tree yields a shape tree. Let   $P_{\mathcal{S}}$ denote the induced probability distribution on the set of shape trees. It
%follows from the same work of \citet{Tajima6628982}. On the other hand, if labels are added to the leaves of a shape tree one obtains phylogenetic trees. Finally, adding leaf labels to ranked trees one obtains labelled ranked trees. The probability distributions $P_{\mathcal{P}}$ on the set of phylogenetic trees and $P_{\mathcal{L}}$ on the set of labelled ranked trees are determined below.

\bigskip
%\textbf{Shape trees:}
%Let $t\in \mathcal{S}$ and let $\gamma(t)$ be the number of internal nodes $i$ such that the sizes of the right and left subtrees of $i$ are different. 
%Moreover, let 
%$$
%\pi(t)=\prod_{i\in J(t)}(n_i-1)\,,
%$$  
%where $n_i$ represents the size of the subtree of internal node $i$ and $J(t)$ is the set of internal nodes of $t$.
%For example, if we consider tree (a) in Fig. \ref{examples}, we have $\gamma=3$ and $\pi=1\times 1\times 2\times 4\times 5 = 40$.
%Given $t \in \mathcal{S}_n$, from \citet{Tajima6628982} follows that 
%$$
%P_{\mathcal{S}}(t)=\frac{2^{\gamma(t)}}{\pi(t)}\,,
%$$
%i.e. the probability of any tree $t \in \mathcal{S}_n$ depends on the parameter $\gamma$ and the
%subtree sizes of all internal nodes, $P_{\mathcal{S}}(t) = P_{\mathcal{S}}(t(\gamma,n_1,...,n_{n-1}))$. 

%\bigskip
\textbf{Probability distribution of ranked trees}

Let $t\in \mathcal{R}$ and let $o(t)$ be the number of internal nodes $i$ whose children are  two leaves. Such internal nodes are called the \emph{cherries} of the tree.
For example, (see Fig.~\ref{ranked6}).
Given $t \in \mathcal{R}_n$, from \citet{Tajima6628982} follows that 
\begin{equation}\label{tagima}
P_{\mathcal{R}}(t)=\frac{2^{n-1-o(t)}}{(n-1)!}\,,
\end{equation}
i.e. the probability of any ranked tree $t \in \mathcal{R}_n$ depends only on two parameters, $o$ and $n$. %$P_{\mathcal{R}}(t) = P_{\mathcal{R}}(t(\alpha,n))$.
\bigskip

\textbf{The probability of generating the same ranked trees twice}

%We consider two cases: first, independently generated trees and, second, trees  linked on a common chromosome. 
Considering trees  linked on a common chromosome one observes  that chromosomal linkage substantially increases the probability that two 'neighboring' trees are identical even if separated by a recombination event. To quantify the effect of linkage and recombination it is important to know the background probability that two independently generated trees are identical.
This probability can be found with the help of the genarating function
%In Fig.~\ref{twopics} we show the result of two simulations of $100$ ranked trees of size $6$  using the program \emph{ms} \citep{hudson:2002}. When two adjacent trees are different they are represented by circles lying on the $x$-axis, when they are identical they are represented by circles above the $x$-axis. The left panel shows the independent case, the right panel shows the case of linkage.
%\begin{figure}
%\begin{center}
%\begin{tabular}{ccc}
%\includegraphics*[angle=0,scale=.33,trim=0 0 0 0]{ind.dat4pic.ps} & ~~~ &
%\includegraphics*[angle=0,scale=.33,trim=0 0 0 0]{rec.new.ps} 
%\end{tabular}
%\end{center}
%\caption{Matching (circles above the line) vs. mismatching (circles on the line) tree pairs. Simulation of $100$ ranked trees of size six. Left: independent case; right: chromosomal linkage.}\label{twopics}
%\end{figure}
%Note that chromosomal linkage substantially increases the probability that two 'neighboring' trees are identical even if separated by a recombination event. To quantify the effect of linkage and recombination it is important to know the background probability that two independently generated trees are identical.
%This probability can be found with the help of the genarating function

$$Y(x,z)=\sum_{t \in \mathcal{R}} \frac{ x^{o(t)} z^{n(t)-1} }{(n(t)-1)!},$$

discussed in more details in Section~\ref{recursive}, eq.~(\ref{Y}).

We have the following result.

\begin{prop}\label{prima}
The probability that two independently generated ranked trees of size $n$ are identical is
\begin{equation}\nonumber
p_n= \frac{4^{n-1}}{(n-1) !} \times [z^{n-1}]Y\left(\frac{1}{4},z\right)\,.
\end{equation}
\end{prop}
\emph{Proof:}
%Let $t\in \mathcal{R}_n$. 
From eq. (\ref{tagima})
%we know that its probability is $P_{\mathcal{R}}(t)=\frac{2^{\alpha(t)}}{n-1 \, !}$.
the probability that $t_1, t_2 \in \mathcal{R}_n$ are identical is

\begin{eqnarray}\nonumber
p_n &=& \sum_{t\in \mathcal{R}_n}P_{\mathcal{R}}(t)^2 \\\nonumber
&=& \frac{1}{(n-1)! ^2}\sum_{t\in \mathcal{R}_n}4^{n-1-o(t)} \\\nonumber
&=& \frac{4^{n-1}}{(n-1)! ^2}\sum_{t\in \mathcal{R}_n}\left(\frac{1}{4}\right)^{o(t)} \\\nonumber
%&=&\frac{1}{\left(\frac{n!(n-1)!}{2^{n-1}}\right)^2}\sum_{t\in R_n}m_2(t)^2 \\\nonumber
%&=&\frac{n!^2}{\left(\frac{n!(n-1)!}{2^{n-1}}\right)^2}\sum_{t\in R_n}\left(\frac{1}{4}\right)^{s(t)} \\\nonumber
%&=&\frac{n!^2(n-1)!}{\left(\frac{n!(n-1)!}{2^{n-1}}\right)^2} \times [z^{n-1}]Y\left(\frac{1}{4},1,1,z\right)\\\nonumber
&=& \frac{4^{n-1}}{(n-1) !} \times [z^{n-1}]Y\left(\frac{1}{4},z\right), \nonumber
\end{eqnarray}

where $[z^{n-1}]Y(1/4,z)$ means the $(n-1)$-st coefficient of the Taylor expansion of $Y(1/4,z)$ in $z=0$. $\Box$

\section{Enumerative results}\label{problems}

\subsection{\small{Outdegree of the  nodes in ranked 
and $0$-$1$-$2$-increasing trees}}

Let $t \in \mathcal{R}_{n}$ and $m=n-1$. Remove all leaves and external branches from $t$ and obtain a reduced tree $\rho(t)$. The tree $\rho(t)$ is a so-called \emph{$0$-$1$-$2$-increasing tree} of size $m$, where, this time, the size is the total number of nodes in the tree and not only of the leaves.
The class $\mathcal{I}_{012}$ of $0$-$1$-$2$-increasing trees is composed of un-ordered rooted trees where all nodes have outdegree $0$, $1$ or $2$. The $m$ nodes of such a tree carry totally ordered labels belonging to $\{1,2,...,m \}$. Moreover, the labelling is such that any child node label is greater than that of the parent node. As usual $\mathcal{I}_{{012}_m}$ denotes the set of $0$-$1$-$2$-increasing trees of size $m$.
Hence, the function $\rho$ is a bijection from $\mathcal{R}_{n}$ to $\mathcal{I}_{{012}_m}$. 

Given a ranked tree $t$, the \emph{outdegree of an internal node} of $t$ is the outdegree of the corresponding node in $\rho(t)$. Thus, if $t \in \mathcal{R}$, the nodes of outdegree $0$ (resp. $1$, $2$) are defined as the nodes with $2$ (resp. $1$, $0$) leaves as direct descendants. 
%In particular the nodes attached to two leaves are usually called the \emph{cherries} of the tree $t$.
%Furthermore, if $t \in \mathcal{R}_{n+1}$, and if $o(t)$ denotes its number of cherries, then 
%$$
%\alpha(t)=n-o(t)\,.
%$$

\bigskip

%According to \cite{MR1251994} the class of $0$-$1$-$2$-increasing trees, and therefore the class of ranked trees via  the bijection $\rho$, is a \emph{variety} of unordered increasing trees characterized by the polynomial 
%$$
%\phi(w)=1+w+\frac{w^2}{2}.
%$$

%In \cite{MR1251994} the authors find several statistics on $0$-$1$-$2$-increasing trees including the one we are going to show in the next section using a different methodology. 

Here, we derive the enumeration of $0$-$1$-$2$-increasing trees with respect to the size and to the number of nodes with outdegree $0$,$1$ and $2$.
The bijection $\rho$ will allow us to use this enumerative result in Section~\ref{weightedF} to determine the probability distribution of the random variable $o$, the \emph{number of cherries}, when $t$ is a ranked tree of size $n$ generated by the coalescent process. It is already known (see \citet{Steel99}) that $o(t)$  is asymptotically normal for large $n$.

\subsubsection{\small{Recursive construction of $0$-$1$-$2$-increasing trees}} \label{recursive}

We show now how the class of $0$-$1$-$2$-increasing trees  can be generated recursively. In particular we construct each tree belonging to $\mathcal{I}_{{012}_{m+1}}$  by adding a new node to some tree in $\mathcal{I}_{{012}_m}$. This construction, denoted by $\Theta$, will then be translated into a functional equation. Solving the equation we obtain a bivariate exponential generating function counting the considered increasing trees with respect to size and to the number of nodes with outdegree $0$, $1$ and $2$.

\begin{figure}
\begin{center}
\includegraphics*[scale=.38,trim=0 0 0 0]{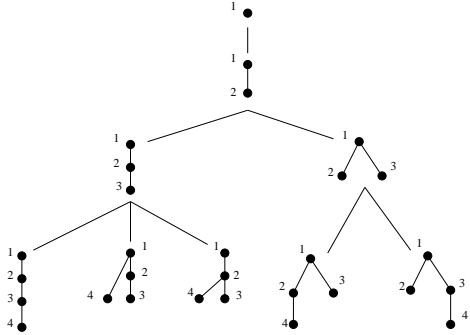}
\end{center}
\caption{First levels of the generating tree associated to $\Theta$.}\label{gentree}
\end{figure}

Given a tree $t \in \mathcal{I}_{{012}_m}$, $\Theta$ simply adds the node labelled '$m+1$' as a child of a node of $t$ having outdegree less than two. Let $o(t),p(t)$ and $q(t)$ denote the number of nodes with outdegree $0,1$ and $2$ respectively. $\Theta$ applied to $t$ produces $o(t)+p(t)$ elements of $\mathcal{I}_{{012}_{m+1}}$ each time adding the new node labelled $m+1$ as a child of the nodes counted in $o(t)+p(t)$.
In Fig.~\ref{gentree} we depict the first steps of this construction process.

Note that $o(t)=q(t)+1$ and $o(t)+p(t)+q(t)=m$. From these relations we have, in particular, that $p(t)=m-2o(t)+1$.
The construction $\Theta$ can be translated into the following \emph{succession rule} (see \citet{BBDFGG}) where each tree is represented by a label composed of the values of its parameters $o$ and $m$ while the exponents show how many times the label is produced, 

$$(o,m) \rightarrow (o,m+1)^o \, (o+1,m+1)^{m-2o+1}.$$

In particular, given a tree $t$ with parameters $o=o(t)$ and $m=m(t)$, the application of $\Theta$ to $t$ produces $o$ new trees having size $m+1$ and $o$ cherries and $m-2o+1$ new trees having size $m+1$ and $o+1$.  
The starting point of the construction is the unique tree of size one represented by $(1,1)$.   

Now consider the exponential generating function $$Y(x,z)=\sum_{t \in \mathcal{I}_{012}} \frac{ x^{o(t)} z^{m(t)} }{m(t) \, !}.$$

The previous succession rule can be translated as follows into an equation for $Y(x,z)$.

\begin{eqnarray} \nonumber
Y(x,z)&=&xz + \sum_{x^o z^m \in \mathcal{I}_{012}} \frac{o x^o z^{m+1}}{(m+1) \, !} + \sum_{x^o z^m \in \mathcal{I}_{012}} \frac{(m-2o+1)(x^{o+1}z^{m+1})}{(m+1) \, !} \\\nonumber
&=& xz + (1-2x)\sum_{x^o z^m \in \mathcal{I}_{012}} \frac{o x^o z^{m+1}}{(m+1) \, !} + xz \sum_{x^o z^m \in \mathcal{I}_{012}} \frac{x^{o}z^{m}}{m \, !} \\\nonumber
&=& xz + (1-2x)\sum_{x^o z^m \in \mathcal{I}_{012}} \frac{o x^o z^{m+1}}{(m+1) \, !} + xz \, Y(x,z) \nonumber
\end{eqnarray}

From the previous equation we obtain that

\begin{eqnarray} \nonumber
\frac{Y(x,z)(1-xz) -xz}{1-2x}&=&\sum_{x^o z^m \in \mathcal{I}_{012}} \frac{o x^o z^{m+1}}{(m+1) \, !}.
\end{eqnarray}

Differentiating both sides with respect to the variable $z$ we have

\begin{eqnarray} \nonumber
\frac{1}{1-2x} \left( \frac{dY}{dz}(x,z) \, (1-xz) -x \, Y(x,z) -x \right)  &=& x \, \frac{dY}{dx}(x,z),
\end{eqnarray}

which is equivalent to

\begin{eqnarray} \label{caso1}
x(1-2x) \, \frac{dY}{dx}(x,z) + (xz-1) \, \frac{dY}{dz}(x,z) &=& -x \, Y(x,z) -x. 
\end{eqnarray}

The previous first order partial differential equation can be solved using the \emph{method of characteristics} (see \cite{CourantHilbert1989}) respecting the condition given by eq.~(\ref{enumR}) $$Y(1,z)=\sec (z)+ \tan (z)-1.$$ Indeed  $Y(1,z)$ must represent the exponential generating function counting $0$-$1$-$2$-increasing trees with respect to size. 

Applying the method consists, first, of solving the two following ordinary differential equations

\begin{eqnarray}\nonumber
z' &=& \frac{xz-1}{x(1-2x)} \\\nonumber
Y' &=& \frac{-xY -x}{x(1-2x)} \nonumber 
\end{eqnarray} 

The solutions are 

\begin{eqnarray} \nonumber
z &=& \frac{c_1 + 2 \, \arctan (\sqrt{2x-1})}{\sqrt{2x-1}}, \\ \label{five}
Y &=& c_2 \, \sqrt{2x-1} -1, \\\nonumber 
\end{eqnarray}

with constants $c_1$ and $c_2$ and where $c_2$ can be written as a function of $c_1$ in the following way

$$c_2=G(c_1)=G( z \, \sqrt{2x-1} -2 \, \arctan (\sqrt{2x-1})).$$

In this way equation (\ref{five}) becomes

$$Y(x,z)=G( z \, \sqrt{2x-1} -2 \, \arctan (\sqrt{2x-1})) \, \sqrt{2x-1}-1,$$

which gives

$$\sec (z)+ \tan (z)-1 \, = \, Y(1,z) \, = \, G(z-\frac{\pi}{2})-1.$$

Function $G$ must satisfy

$$G(z)= \sec (z+\frac{\pi}{2})+ \tan (z+\frac{\pi}{2})=\frac{-1- \cos (z)}{\sin (z)}.$$

Inserting this into (\ref{five}) we have 

$$Y(x,z)=\sqrt{2x-1} \left( \frac{-1- \cos (z \, \sqrt{2x-1} -2 \, \arctan (\sqrt{2x-1}))}{\sin (z \, \sqrt{2x-1} -2 \, \arctan (\sqrt{2x-1}))} \right) -1,$$

which, after some calculations, finally gives 

\begin{equation}\label{ipsilon1}
Y(x,z)=\frac{\sqrt{2x-1}}{\tan \left( -\frac{z\sqrt{2x-1}}{2}+ \arctan (\sqrt{2x-1}) \right) }  - 1.
\end{equation}

%The function $Y$ has a Taylor expansion at $z=0$ whose first terms are:

%\begin{eqnarray} \nonumber
%Y(x,z) &=& xz \\\nonumber
%&& + \frac{1}{2} xz^2 \\\nonumber
%&& + \frac{1}{6} (x^2 + x)z^3 \\\nonumber
%&& + \frac{1}{24} (4 x^2 + x)z^4 \\\nonumber
%&& + \frac{1}{120} (4 x^3 + 11 x^2 + x)z^5 \\\nonumber
%&& + \frac{1}{720} (34 x^3 + 26 x^2 + x)z^6 \\\nonumber
%&& + \frac{1}{5040} (34 x^4 + 180 x^3 + 57 x^2 + x)z^7 \\\nonumber
%&& + ... \, .   \nonumber
%\end{eqnarray}

Note that the condition $Y(1,z)=\sec (z)+ \tan (z)-1$ is respected. 

Indeed 

\begin{eqnarray}\nonumber
Y(1,z)= \frac{1}{\tan \left( -\frac{z}{2} + \frac{\pi}{4} \right) } -1
\end{eqnarray}
 
and

\begin{eqnarray} \nonumber
\frac{1}{\tan \left( -\frac{z}{2} + \frac{\pi}{4} \right) } 
&=& \frac{1+ \tan \left( \frac{z}{2} \right) }{1- \tan \left( \frac{z}{2} \right)} 
= \frac{1+ \cos (z)+ \sin (z)}{1+ \cos (z)- \sin (z)} \\\nonumber 
&=&\frac{1+{\cos}^2(z)+2 \cos (z)-{\sin}^2(z)}{(1+ \cos (z)- \sin (z))^2} \\\nonumber 
&=& \frac{\cos (z)}{1- \sin (z)} 
= \frac{1+ \sin (z)}{\cos (z)}
\end{eqnarray}

Moreover, using the fact that $$\exp{(z \sqrt{-2x+1})}=\cos (z \sqrt{2x-1})+ i \, \sin (z \sqrt{2x-1}),$$ we can write eq.~(\ref{ipsilon1})  in terms of the exponential function as

\begin{equation} \label{Y}
Y(x,z)= \frac{2 \, { \left( x \exp{ \left( \sqrt{-2 \, x + 1} z \right) } - x \right) }}{{ \left( \sqrt{-2 \, x + 1} - 1 \right) } \exp{ \left( \sqrt{-2 \, x + 1} z \right) }  + \sqrt{-2 \, x + 1} + 1}.
\end{equation}

Performing the substitution $x=1/4$ we have that  

$$
Y\left(\frac{1}{4},z\right)=\frac{e^{\left(\sqrt{\frac{1}{2}} z\right)} - 1}{2 \,
{\left({\left(\sqrt{\frac{1}{2}} - 1\right)} e^{\left(\sqrt{\frac{1}{2}}
z\right)} + \sqrt{\frac{1}{2}} + 1\right)}}\, ,
$$

the Taylor expansion of which is 

$$
Y\left(\frac{1}{4},z\right)=\frac{1}{4}z+\frac{1}{8}z^2+\frac{5}{96}z^3+\frac{1}{48}z^4+\frac{1}{120}z^5+\dots.$$
Using the result of Proposition~\ref{prima} we can now effectively calculate the probability $p_n$ that two ranked trees having $n$ leaves are identical when generated independently by the coalescent process:
$p_2=\frac{4}{1!} \times \frac{1}{4}=1, p_3=\frac{4^2}{2!}  \times \frac{1}{8}=1,p_4=\frac{4^3}{3!} \times \frac{5}{96}=\frac{5}{9},p_5=\frac{4^4}{4!}  \times \frac{1}{48}=\frac{2}{9}$ and $p_6=\frac{4^5}{5!} \times  \frac{1}{120}=\frac{16}{225}$, and so on. 

%In Fig.~\ref{twopics} we show the result of two simulations of $100$ ranked trees of size $6$  using the program \emph{ms} \citep{hudson:2002}. We considered two cases: first, independently generated trees and, second, trees  linked on a common chromosome. When two subsequent trees are different they are represented by circles lying on the $x$-axis, when they are identical they are represented by circles above the $x$-axis. The left panel shows the indepenedent case, the right panel shows the case of linkage.

%\begin{figure}
%\begin{center}
%\begin{tabular}{ccc}
%\includegraphics*[angle=0,scale=.33,trim=0 0 0 0]{ind.dat4pic.ps} & ~~~ &
%\includegraphics*[angle=0,scale=.33,trim=0 0 0 0]{rec.dat4pic.ps} 
%\end{tabular}
%\end{center}
%\caption{Matching (circles above the line) vs. mismatching (circles on the line) tree pairs. Simulation of $100$ ranked trees of size six. Left: independent case; right: chromosomal linkage. For further explanation see the text.}\label{twopics}
%\end{figure}

%Note that chromosomal linkage substantially increases the probability that two 'neighboring' trees are identical. To quantify the effect of linkage it is important to know the background probability that two independently generated trees are identical.

\subsubsection{\small{The probability distribution of the number of cherries}} \label{weightedF}

We are now ready to state the enumeration of ranked trees with respect to size and number of nodes of outdegree $0$, $1$ or $2$, when each tree is weighted by its probability under the coalescent process. This exact enumerative result is novel and  
achieved with the help of the weighted generating function
\begin{equation} \nonumber
F(x,z)=\sum_{t\in \mathcal{R}_n, \, n>1} \frac{2^{n(t)-1-o(t)}}{(n(t)-1)!}x^{o(t)}z^{n(t)}.
\end{equation}  

Function $F$ has a more intuitive interpretation if one considers the transformation $Y_w=\frac{F}{z}$ instead. 
It can be interpreted as a weighted exponential generating function counting $0$-$1$-$2$-increasing trees with respect to the outdegree and the total number of nodes.

Starting from equation (\ref{Y}), we perform some substitutions on  $Y$ to obtain $Y_w$. 
%In particular, we introduce two new variables $v$ and $y$ and define $\tilde{Y}
%(x,v,z,y)=Y(xv,zy)$. 
In particular we have $Y_w = Y \left( \frac{x}{2},2z \right) $ and, multiplying by $z$, we finally obtain the desired function $F$.

\begin{prop}\label{pere}
The weighted ordinary generating function of ranked trees considered with respect to size and number of cherries is 
\begin{equation} \label{oggi}
F(x,z)= \frac{zx \exp{\left(2z \, \sqrt{-x + 1}\right)} - zx}{{\left( \sqrt{-x  + 1} -1\right)} \exp{\left(2z \, \sqrt{-x  + 1}\right)} + 1 + \sqrt{-x  + 1}}.
\end{equation}
The probability of having $o'$ cherries in a ranked tree of size $n$ corresponds to the coefficient of $x^{o'}z^n$ in the Taylor expansion of $F$ around $z=0$, i.e. $$P_n(o= o')=[x^{o'}z^n]F(x,z).$$
\end{prop}

%The probability of having $\overline{o}$ cherries in a ranked tree of size $n$ corresponds to the coefficient of $x^{\overline{o}}z^n$ in the Taylor expansion of $F$ centered at $z=0$, i.e. $$P_n(o=\overline{o})=[x^{\overline{o}}z^n]F(x,z).$$ Extracting the coefficients we have the next result. 

%\begin{coro}\label{reccur}
%The following recursion holds:
%\begin{itemize}
%\item[i)] If  $\overline{o}=0$ or $\overline{o} > \lfloor n/2 \rfloor$
%\begin{equation} \nonumber
%P_n(o=\overline{o})  =  0. 
%\end{equation}
%\item[ii)] $P_2(o=1)=1$ while, if $0<\overline{o} \leq  \lfloor n/2 \rfloor$, we have
%\begin{equation} \nonumber
%P_n(o=\overline{o}) = \frac{2\overline{o}}{n-1} P_{n-1}(o=\overline{o}) + \frac{n+1-2\overline{o}}{n-1} P_{n-1}(o=\overline{o}-1).
%\end{equation}
%\end{itemize}
%\end{coro}

The first terms of the Taylor expansion of (\ref{oggi}) are described below;
  
\begin{eqnarray} \nonumber
F(x,z)&=&x z^2 \\\nonumber
&&+ x z^{3} \\\nonumber
&&+ \frac{1}{3} \, {\left(x^{2} + 2 \, x \right)} z^{4} \\\nonumber
&&+ \frac{1}{3} \, {\left(2 \, x^{2} + x \right)} z^{5} \\\nonumber 
&&+ \frac{1}{15} \, {\left(2 \,
x^{3}  + 11 \, x^{2} + 2 \, x\right)} z^{6} \\\nonumber
&&+ \frac{1}{45} \,
{\left(17 \, x^{3}  + 26 \, x^{2}
+ 2 \, x \right)} z^{7} \\\nonumber
&&+ \frac{1}{315} \, {\left(17 \, x^{4}
 + 180 \, x^{3} + 114 \, x^{2}
 + 4 \, x \right)} z^{8}  \\\nonumber
&& + \dots \, . \nonumber
\end{eqnarray}    

Looking at Fig.~\ref{ranked6} one can check that, for example, there are exactly $11$ trees represented by the monomial $x^2z^6$. Each one of them has probability $\frac{1}{15}$. This is in agreement with the term $\frac{11}{15}x^2z^6$ in the expansion. Indeed, $\frac{11}{15}$ is the probability to obtain a ranked tree of size $6$ with two cherries.

%\begin{figure}
%\begin{center}
%\includegraphics*[scale=.4,trim=0 0 0 0]{alberi.eps}
%\end{center}
%\caption{The sixteen ranked trees of size $6$. Each row is an equivalence class of exactly one shape.}\label{ranked6}
%\end{figure}

%There are exactly $11$ trees represented by the monomial $x^2z^6$. These are the ones corresponding to the second, third, foth and fifth row of Fig.~\ref{ranked6}. Each one of them has probability $\frac{1}{15}$. This is in agreement with the term $\frac{11}{15}x^2z^6$ in the Taylor expansion. Indeed, $\frac{11}{15}$ is the probability to obtain a ranked tree of size $6$ with two cherries.  

%Analogously, the term $\frac{2}{15}x z^6$ in the expansion reflects the fact that there is exactly one tree corresponding to the monomial $x z^6$ and that this tree has probability $\frac{2}{15}$. 

%Finally, the term $\frac{2}{15}x^3z^6$ is the monomial $x^3z^6$ weighted with the factor $4/30$ and
%corresponds to the exactly four trees with three cherries, each one of them having probability $\frac{1}{30}$.

\bigskip

Using the result of Proposition~\ref{pere} we compute the discrete probability distribution of the random variable $o(t)$ for trees of fixed size $n$. In this case $o$ is a random variable which takes values between $1$ and $\lfloor n/2 \rfloor$. In Fig.~\ref{distpitchfork} we have depicted the distribution of $o$ for a ranked tree of size $n=54$.

%\begin{figure}
%\begin{center}
%\includegraphics*[scale=.5,trim=0 0 0 0]{pofsC_100b.eps}
%\includegraphics*[scale=.4,trim=0 0 0 0]{data4s-plot.eps}
%\end{center}
%\caption{Probability distribution of the random variable $o$ in the case $\mathcal{R}_{100}$.}\label{u0}
%\end{figure}

By Proposition~\ref{pere} one can also determine the expected value $E_o(n)$ and the variance $Var_o(n)$ of the random variable $o$ in dependence of tree size $n$. Using other methods these have been determined before, for example by \citet{Steel99}.

Using our approach the expectation is
$$
E_o(n)=[z^n]\frac{dF}{dx}(1,z)=[z^n]\frac{z^4-3z^3+3z^2}{3(z-1)^2}\,.
$$
If $n>2$, this simplifies to 
$$
E_o(n)=\frac{n}{3}\,.
$$

%\bigskip 

%In Fig.~\ref{simulation} we have depicted $o$ for $300$ independently generated ranked trees using the program \emph{ms} \citep{hudson:2002}. Each tree is of size $n=21$ and is generated according to the random coalescent process. The $y$-axis shows the realization of $o$ for each of the $300$ trees. The theoretical expectation is $E_o(21)=\frac{21}{3}=7$, which agrees with the computer simulations. 

%\begin{figure}
%\begin{center}
%\includegraphics*[scale=.4]{output.eps}
%\end{center}
%\caption{Simulations for $300$ ranked trees of size $21$. The $y$-axis shows the random variable $o$.}\label{simulation}
%\end{figure}

The second moment is 
\begin{eqnarray} \nonumber
E_{o^2}(n)&=&[z^n]\frac{d(x\frac{dF}{dx})}{dx}(1,z)=[z^n]\frac{d^2F}{dx}(1,z)+ [z^n]\frac{dF}{dx}(1,z) \\\nonumber
&=&[z^n]\frac{2(z^7-6z^6+15z^5-15z^4)}{45(z-1)^3}+E_o(n) \\\nonumber
&=&[z^n]\left(\frac{2}{(z-1)^3}\left(\frac{z^7}{45}-\frac{2z^6}{15}+\frac{z^5}{3}-\frac{z^4}{3}\right)\right) + E_o(n)\,.
\end{eqnarray}
If $n>6$, and using $Var_o(n)=E_{o^2}(n)-E_o^2(n)$, we obtain the variance of $o$
\begin{eqnarray} \nonumber
{\mathrm Var}_o(n)&=&-\frac{(n-5)(n-6)}{45} +\frac{2(n-4)(n-5)}{15} \\\nonumber
&&-\frac{(n-3)(n-4)}{3} +\frac{(n-2)(n-3)}{3} \\\nonumber
&&+\frac{n}{3} -\frac{n^2}{9} \\\nonumber 
&=& \frac{2n}{45}\,. \nonumber
\end{eqnarray}

\bigskip

%Expecation and Variance of $\alpha$ (see Section~\ref{process}) can be easily derived from the expectation and variance of $o$.

%\begin{prop} \label{E}
%For ranked trees of size $n > 6$ the expected value of the random variable $\alpha$ is 
%$$
%E_{\alpha}(n)=\frac{2n-3}{3},$$ 
%and the variance is 
%$$
%{\mathrm Var}_{\alpha}(n)=\frac{2n}{45}.
%$$
%\end{prop}

%Proposition \ref{E} gives an (partial) answer to our initial question: how does a typical ranked tree look like?
%The variance of the trees shown in Fig~\ref{simulation} is approximately $0.94$, very close to the theoretical value.
Note that this is the variance of cherries of independently generated trees. Considering 'linked' trees, i.e. along a recombining chromosome, the variance is smaller.

\subsection{The number of pitchforks}   

%Trees are very popular combinatorial objects. They have been enumerated with respect to various parameters. In biological applications it is important to consider that trees may be generated by particular random processes, for instance a pure birth process or a coalescent process. In these applications trees are usually not uniformly distributed. 
%In the previous sections we have derived several statistics concerning ranked trees and the random variable \emph{number of cherries} when generated under the coalescent process. It will be of interest to extend 

The recursive construction presented in Section \ref{recursive} can be extended in order to consider also \emph{pitchforks}.
 
%to other parameters and to see whether the resulting functional equation remains solvable. 

Using different methods, they have been studied before for example by \citet{rosenberg}.
%As an example one can consider the number of \emph{pitchforks} which have been studied for example by \citet{rosenberg}. 
A pitchfork in a ranked (resp. $0$-$1$-$2$-increasing) tree is simply a subtree  with $3$ leaves (resp. $2$ nodes).
If $r(t)$ denotes the number of pitchforks in $t \in \mathcal{I}_{012}$ the construction of Section \ref{recursive} is extended  to the new random variable $r$. We find the following succession rule: 
\begin{eqnarray} \nonumber
(o,r,m) &\rightarrow & (o,r,m+1)^r \, (o,r+1,m+1)^{o-r} \\\nonumber 
(o,r,m) &\rightarrow & (o+1,r-1,m+1)^{r} \, (o+1,r,m+1)^{m-2o+1-r}. \nonumber
\end{eqnarray}

Considering now 

$$Y(x,v,z)=\sum_{t \in \mathcal{I}_{012}} \frac{ x^{o(t)} v^{r(t)} z^{m(t)} }{m(t) \, !},$$

we obtain the following differential equation:

\begin{equation}\label{nosol}
(v+x)(v-1)\frac{dY}{dv} = x+ xY + x(v-2x)\frac{dY}{dx} + (xz-1)\frac{dY}{dz}.
\end{equation}

For $v=1$ it reduces to eq.~(\ref{caso1}) but there is non easy analytic solution.
%which gives eq.~(\ref{caso1}) in the case $v=1$ and seems not so easy to solve.

However, we can still obtain the expected value $E_r(m)$ for the number of pitchforks in $0$-$1$-$2$ increasing trees with  $m$ nodes. Starting from (\ref{nosol}) and performing the substitutions $x=1/2$ and $z=2z$ we obtain
%\begin{eqnarray} \nonumber
%\left( v+\frac{1}{2} \right) (v-1)\frac{dY}{dv} \left( \frac{1}{2},v,2z \right) &=& \frac{1}{2}+ \frac{1}{2} \, Y\left( \frac{1}{2},v,2z \right) \\\nonumber
%&& + \frac{1}{2}(v-1)\frac{dY}{dx} \left( \frac{1}{2},v,2z \right) \\\nonumber 
%&& + (z-1)\frac{dY}{dz} \left( \frac{1}{2},v,2z \right), \\\nonumber
%\end{eqnarray}

%from which we have

\begin{eqnarray} \nonumber
\frac{dY}{dv} \left( \frac{1}{2},v,2z \right) &=& \frac{1+Y\left( \frac{1}{2},v,2z \right)+2(z-1)\frac{dY}{dz} \left( \frac{1}{2},v,2z \right)}{2\left( v+\frac{1}{2} \right) (v-1)} \\\nonumber
&& +\frac{\frac{dY}{dx} \left( \frac{1}{2},v,2z \right)}{2\left( v+\frac{1}{2} \right)} \nonumber
\end{eqnarray}

from which we have
\begin{eqnarray} \nonumber
[z^m]\frac{dY}{dv} \left( \frac{1}{2},v,2z \right)&=&[z^m]\frac{Y\left( \frac{1}{2},v,2z \right)+2(z-1)\frac{dY}{dz} \left( \frac{1}{2},v,2z \right)}{2\left( v+\frac{1}{2} \right) (v-1)} \\\nonumber
&& +[z^m]\frac{\frac{dY}{dx} \left( \frac{1}{2},v,2z \right)}{2\left( v+\frac{1}{2} \right)}. \nonumber
\end{eqnarray}

When $v\rightarrow 1$ we find that

\begin{eqnarray} \nonumber
E_r(m) &=& [z^m]\left( \lim_{v\rightarrow 1}\frac{Y\left( \frac{1}{2},v,2z \right)+2(z-1)\frac{dY}{dz} \left( \frac{1}{2},v,2z \right)}{2\left( v+\frac{1}{2} \right) (v-1)}\right) \\\nonumber
&& + \frac{2 E_o(m)}{3}. \nonumber
\end{eqnarray}

The considered limit can be determined according to \emph{l' Hospital's rule} taking the derivative of the numerator and the denominator with respect to $v$ and performing then the substitution $v=1$. Furthermore, from Section~\ref{weightedF} $E_o(m)=(m+1)/3$, and thus

\begin{eqnarray} \nonumber 
E_r(m) &=& [z^m] \left(\frac{1}{3} \sum_{t\in \mathcal{I}_{012}} r(t)\frac{2^{m(t)-o(t)}}{m(t) \, !} z^{m(t)} \right) \\\nonumber
&& + [z^m] \left( \frac{2}{3}(z-1) \sum_{t\in \mathcal{I}_{012}} r(t)m(t)\frac{2^{m(t)-1-o(t)}}{m(t) \, !} z^{m(t)-1} \right) \\\nonumber 
&& + \frac{2(m+1)}{9} \\\nonumber
&=& \frac{1}{3} E_r(m) + [z^m] \left( \frac{z-1}{3z} \sum_{k>0} k E_r(k) z^k \right) +  \frac{2(m+1)}{9}   \\\nonumber
&=& \frac{1}{3} E_r(m) + \frac{mE_r(m)-(m+1)E_r(m+1)}{3} +  \frac{2(m+1)}{9}.   \\\nonumber
\end{eqnarray}

Reordering terms we obtain the recursion

\begin{eqnarray}\nonumber
E_r(2)&=&1 ; \\\nonumber
(m+1)E_r(m+1)&=&(m-2)E_r(m)+\frac{2(m+1)}{3}.\nonumber
\end{eqnarray}

This gives for an increasing tree with $m>2$ nodes  $$E_r(m)=\frac{m+1}{6}.$$

From eq. (\ref{nosol}) one can also compute the full probability distribution of the random variable $r$ when an increasing tree of fixed size is generated by the coalescent process. Indeed, if we consider 

$$Y_m(x,v,z)= \sum_{t \in  \mathcal{I}_{{012}_m} } \frac{ x^{o(t)} v^{r(t)} z^{m} }{m \, !}$$ 

the following result provides a recursion which can be used to compute the functions $Y_m$ for any $m \geq 1$. 

\begin{prop} 
The following recursion holds:
\begin{small}
\begin{eqnarray} \nonumber
Y_1 &=& xz \\\nonumber
Y_{m+1} &=& \int \left[ (v+x)(1-v)\frac{dY_m}{dv} + xY_m + x(v-2x)\frac{dY_m}{dx} + xz \frac{dY_m}{dz} \right] dz \nonumber
\end{eqnarray}
\end{small} 
\end{prop}
\emph{Proof.} Consider eq.~(\ref{nosol}) without the monomial $x$ which appears there. If we then isolate the term  $\frac{dY}{dz}$ and integrate both sides of the resulting equation with respect to the variable $z$ we obtain the polynomial $Y_{m+1}$ starting from $Y=Y_m$. $\Box$

\medskip 

The results for $m=1,2,3,4,5$ are as follows

\begin{eqnarray}\nonumber
Y_1 &=& xz \\\nonumber
Y_2 &=& \frac{1}{2} v x z^2 \\\nonumber
Y_3 &=& \frac {1} {6} v x z^3 + \frac {x^2 z^3} {6} \\\nonumber
Y_4 &=& \frac{1}{24} v x z^4 + \frac{x^2 z^4}{24} + \frac{1}{8} v x^2 z^4 \\\nonumber
Y_5 &=& \frac{1}{120} v x z^5 + \frac{x^2 z^5}{120} + \frac{7}{120} v x^2 z^5 + 
\frac{1}{40} v^2 x^2 z^5 + \frac{x^3 z^5}{30} \nonumber
\end{eqnarray}

\begin{figure}
\begin{center}
\includegraphics*[scale=.4,trim=0 0 0 0]{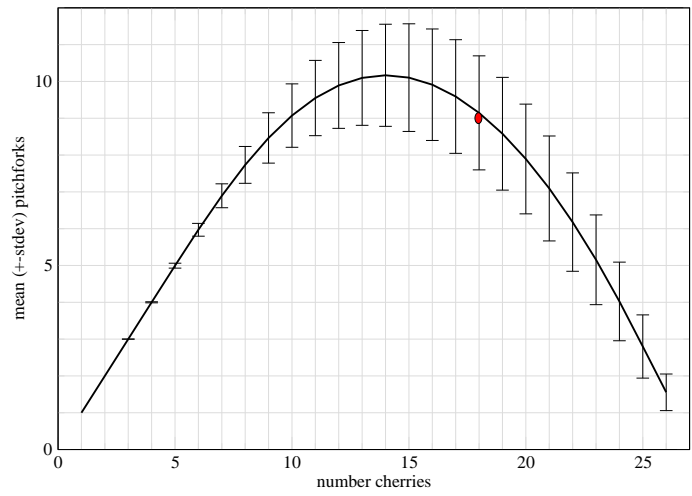}
\end{center}
\caption{Mean of the conditional probability distribution of pitchforks given the number of cherries for  $\mathcal{R}_{54}$.}\label{joint}
\end{figure}

The above results concerning cherries and pitchforks can be extended  to the joint and conditional probability distributions (see Fig.~\ref{joint}). Summarizing, we state

\begin{figure}
\begin{center}
\includegraphics*[scale=.4,trim=0 0 0 0]{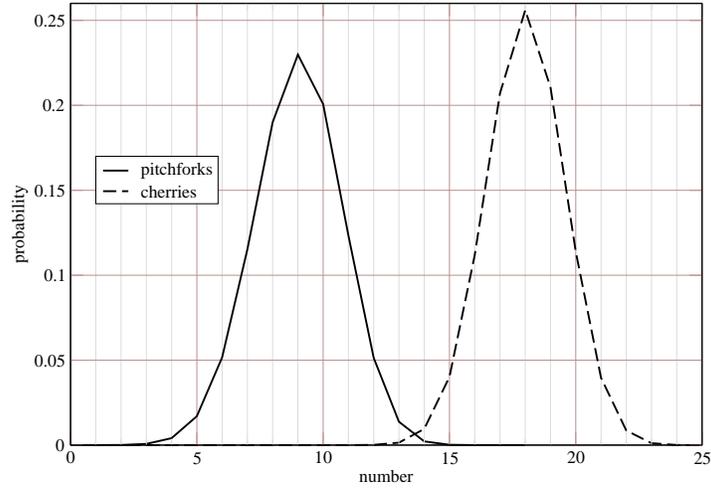}
\end{center}
\caption{ Distributions of cherries and pitchforks for $\mathcal{R}_{54}$ (i.e. $\mathcal{I}_{{012}_{53}}$).}\label{distpitchfork}
\end{figure}

\begin{prop} 
\begin{itemize}
\item[i)]The probability of having $r'$ pitchforks in an increasing tree of size $m$ (see Fig.~\ref{distpitchfork}) is 
$$P_m(r= r') = [v^{r'}]Y_m \left( \frac{1}{2},v,2 \right);$$
\item[ii)] The probability of having $o'$ cherries and $r'$ pitchforks in an increasing tree of size $m$ is 
$$P_m(o= o',r= r') = [x^{o'}v^{r'}]Y_m \left( \frac{x}{2},v,2 \right);$$
\item[iii)] The probability of having $r'$ pitchforks in an increasing tree of size $m$ given it has $o'$ cherries (see  Fig.~\ref{joint}) is 
$$P_m(r= r' | o= o') = \frac{P_m(o= o',r= r')}{P_m(o= o')} = \frac{[x^{o'}v^{r'}]Y_m \left( \frac{x}{2},v,2 \right)}{ [x^{o'}]Y_m \left( \frac{x}{2},1,2 \right)}.$$
\end{itemize}
\end{prop}

%Furthermore, it will be of interest to see whether similar properties can be derived
%for shape and phylogenetic trees. 
%A promising concept in this context is that of tree isomorphism. Two trees are isomorphic
%if they have identical shapes. Some results along these lines have recently been obtained by \citet{Bona:1154505}.
%In a population genetic framework it will be interesting to further explore the effect of recombination upon the distributions of tree shape statistics. 

%\begin{figure}
%\begin{center}
%\includegraphics*[scale=0.115,angle=0]{pezzetti.eps}
%\end{center}
%\caption{Cardinalities of $\mathcal{S}_5,\mathcal{R}_5,\mathcal{P}_5$ and $\mathcal{L}_5$.}\label{sizes5}
%\end{figure}

\subsection*{Acknowledgments}
We gratefully acknowledge helpful discussions with L. Ferretti, A. Klassmann and A. Malina. Financial support was provided by the German Research Foundation (DFG-SFB680).

\bibliographystyle{plain}

\end{document}